\documentclass[12pt]{article}
\title{\bf Infinite divisibility of Smith matrices
\author{{\sc Shaofang Hong} (Chengdu)}
\thanks{{\it Key words and Phrases:} multiplicative function, $l$-th Dirichlet convolution,
Smith matrix, infinite divisibility.}%
\thanks{{\it 2000 Mathematics Subject Classification:} Primary 11C20, 11A25}
\thanks{The research was supported by Program for New Century Excellent Talents
in University Grant \# NCET-06-0785}}%
\date{}

\begin{document}
\baselineskip=17pt
\maketitle
\newcommand{\SSS}{\stackrel}
\newcommand{\p}[2]{{\Phi_{#1}}(#2) }
\newcommand{\J}[2]{\left( \frac{#1}{#2}\right)}
\renewcommand{\thesection}{\arabic{section}}
{\sc Abstract.} Given an arithmetical function $f$, by $f(a, b)$ and
$f[a, b]$ we denote the function $f$ evaluated at the greatest
common divisor $(a, b)$ of positive integers $a$ and $b$ and
evaluated at the least common multiple $[a, b]$ respectively. A
positive semi-definite matrix $A=(a_{ij})$ with $a_{ij}\ge 0$ for
all $i$ and $j$ is called infinitely divisible if the fractional
Hadamard power $A^{\circ r}=(a_{ij}^r)$ is positive semi-definite
for every nonnegative real number $r$. Let $S=\{x_1, ..., x_n\}$ be
a set of $n$ distinct positive integers. In this paper, we show that
if $f$ is a multiplicative function such that $(f*\mu )(d)\ge 0$
whenever $d|x$ for any $x\in S$, then the $n\times n$ matrices
$(f(x_i, x_j))$, $(\frac{1}{f[x_i, x_j]})$ and $(\frac{f(x_i,
x_j)}{f[x_i, x_j]})$ are infinitely divisible. Finally we extend
these results to the Dirichlet convolution case which produces
infinitely many examples of infinitely divisible matrices. Our
results extend the results
obtained previously by Bourque, Ligh, Bhatia, Hong, Lee, Lindqvist and Seip.\\

{\bf 1. Introduction.} Given an arithmetical function $f$, by $f(a,
b)$ and $f[a, b]$ we denote the function $f$ evaluated at the
greatest common divisor $(a, b)$ of positive integers $a$ and $b$
and evaluated at the least common multiple $[a, b]$ respectively. In
1875, Smith [18] showed his renowned result stating that the
determinant of the $n\times n$ matrix $[f(i,j)]$, which has $f(i,j)$
as its $(i,j)$-entry, is the product $\prod^n_{k=1}(f*\mu )(k)$,
where $\mu $ is M\"obius function and $f*\mu $ is the Dirichlet
convolution of $f$ and $\mu $ defined for any integer $a\ge 1$ by
$(f*\mu )(a)=\sum_{d|a}f(d)\mu (a/d)$, where $d$ runs over all
positive divisors of $a$. Since then many generalizations and
related results have been published. See, for instance, [1, 5-17,
19]. Later on, all such kind of matrices are called {\it Smith
matrices}.

A positive semi-definite matrix $A=(a_{ij})$ with $a_{ij}\ge 0$ for
all $i$ and $j$ is called {\it infinitely divisible} if the {\it
fractional Hadamard power} $A^{\circ r}=(a_{ij}^r)$ is positive
semi-definite for every nonnegative real number $r$. Infinitely
divisible matrices arise in several different contexts. Bhatia [3]
and Bhatia and Kosaki [4] discussed this topic and presented some
examples of infinitely divisible matrices. Throughout this paper we
always let $S=\{x_1, ..., x_n\}$ be any given set of $n$ distinct
positive integers. Bourque and Ligh [5] showed that the power GCD
matrix $((x_i, x_j)^r)$ defined on $S$ is positive definite if
$r>0$. From this one can read immediately that the {\it power GCD
matrix} $((x_i, x_j)^r)$ is infinitely divisible if $r>0$. It
follows from [13] that the {\it reciprocal power LCM matrix}
$(\frac{1}{[x_i, x_j]^r})$ is positive definite, and hence is
infinitely divisible if $r>0$. Note that an {\it LCM matrix} $([x_i,
x_j])$ and a {\it reciprocal GCD matrix} $(\frac{1}{[x_i, x_j]})$
may be singular (see [9, 11]). Bourque and Ligh [5] proved that the
matrix $(f(x_i, x_j))$ is positive definite if $f\in \tilde{\cal
C}_S:=\{f|(f*\mu )(d)>0 \ {\rm whenever} \ d|x \ {\rm for \ any} \
x\in S\}$.

In this paper, we consider a natural class of arithmetical functions
${\cal C}_S:=\{f|(f*\mu )(d)\ge 0 \ {\rm whenever} \ d|x \ {\rm for
\ any} \ x\in S\}$. Using the continuity argument we show that the
matrix $(f(x_i, x_j))$ is positive semi-definite if $f\in {\cal
C}_S$. Consequently we show that such matrix is infinitely divisible
if $f\in {\cal C}_S$ is multiplicative. But it fails to be true if
$f$ is not multiplicative. We show also that the matrices
$(\frac{1}{f[x_i, x_j]})$ and $(\frac{f(x_i, x_j)}{f[x_i, x_j]})$
are infinitely divisible if $f\in {\cal C}_S$ is multiplicative,
where we make the convention $1/f(a):=0$ if $f(a)=0$. Finally we
extend these results to the Dirichlet convolution case which
produces infinitely many examples of infinitely divisible matrices.
Our results extend the results of Bourque, Ligh, Bhatia, Hong, Lee,
Lindqvist and Seip. We refer the readers to [2] for the basic
elementary concepts and facts from number theory.\\

{\bf 2. Lemmas and theorems.} By [8] and using the continuity
argument we show the following result.\\

{\sc Lemma 2.1.} {\it If $f\in {\cal C}_{S}$, then the $n\times n$
matrix $((f(x_i, x_j))$ is positive semi-definite.}

{\it Proof.} Let $f\in {\cal C}_{S}$ and pick $\epsilon >0$ and
$\bar f\in \tilde{\cal C}_{S}$. Then it is easy to see that
$f+\epsilon \bar f\in \tilde{\cal C}_{S}$. For an arithmetical
function $g$ and $1\le k\le n$, let
$$
\alpha _g(x_k):=\sum_{\SSS{d|x_k} {d\dagger x_t, \ x_t<x_k}}(g*\mu
)(d).
$$
By
Theorem 1 of [8],
$${\rm det}((f+\epsilon \bar f)(x_i, x_j))\ge \displaystyle
\prod_{k=1}^n\alpha _{f+\epsilon \bar f}(x_k). \ \ \ \eqno (1)$$

Note that both sides of (1) are polynomials in $\varepsilon $.
Moreover, the constant coefficients of the left and right hand sides
are, respectively, ${\rm det}(f(x_i, x_j))$ and $\prod_{k=1}^n\alpha
_f(x_k)$. Since (1) holds for any $\varepsilon>0$, letting
$\varepsilon\rightarrow 0$ we obtain

$${\rm det}(f(x_i, x_j))\ge \displaystyle
\prod_{k=1}^n\alpha _{f}(x_k). \ \ \ \eqno (2)$$

For any $1\le l\le n$, since $f\in {\cal C}_{S}$, then the
inequality (2) implies that the determinant of any principal
submatrix of order $l$ of $(f(x_i, x_j))$ is nonnegative. This
concludes that the matrix $((f(x_i, x_j))$ is positive
semi-definite. \\

Associated to any nonnegative real number $r$ and the arithmetical
function $f$ such that $f(m)\ge 0$ for any positive integer $m$, we
define the arithmetical function $f^r$ by $f^r(m):=f(m)^r$ for any
positive integer $m$.\\

{\sc  Lemma 2.2.} {\it Let $f$ be a multiplicative function such
that $f\in {\cal C}_S$. Then for any nonnegative real number $r$, we
have $f^r\in {\cal C}_S$.}

{\it Proof.} Since $f$ is multiplicative, so is $f^r$. Hence
$f^r*\mu $ is multiplicative. Therefore it is sufficient to show
that $(f^r*\mu )(p^e)\ge 0$ for any prime $p$ and any positive
integer $e$ such that $p^e$ divides some $x\in S$. This will be done
in the following.

Evidently one has
$$
(f*\mu ) (p^e)(f^r*\mu )(p^e)=(f(p^{e})-f(p^{e-1}))
(f(p^{e})^r-f(p^{e-1})^r). \eqno (3)
$$
Notice that $f(p^e)=\sum_{t=0}^e(f*\mu )(p^t)\ge 0$ for any
nonnegative integer $e$ since $f\in {\cal C}_S$ implying that
$(f*\mu )(p^t)\ge 0$ for any integer $0\le t\le e$. Since $r\ge 0$,
we have $f(p^{e})^r\ge f(p^{e-1})^r\ge 0$ (resp. $0\le f(p^{e})^r\le
f(p^{e-1})^r$) if $f(p^{e})\ge f(p^{e-1})$ (resp. $f(p^{e})\le
f(p^{e-1})$). Thus
$$
(f(p^{e})-f(p^{e-1}))(f(p^{e})^r-f(p^{e-1})^r)\ge 0. \eqno (4)
$$
But $(f*\mu ) (p^e)\ge 0$ since $f\in {\cal C}_S$. Then by (3) and
(4), $(f^r*\mu )(p^e)\ge 0$. So $f^r\in {\cal C}_S$ as required. \\

{\sc Theorem 2.1.} {\it If $f$ is multiplicative and $f\in {\cal
C}_{S}$, then the $n\times n$ matrices $(f(x_i, x_j))$,
$(\frac{1}{f[x_i, x_j]})$ and $(\frac{f(x_i, x_j)}{f[x_i, x_j]})$
are infinitely divisible.}

{\it Proof.} Since $f\in {\cal C}_{S}$ is multiplicative, by Lemma
2.2 we know that for any $r\ge 0$, $f^r\in {\cal C}_{S}$. It then
follows from Lemma 2.1 that for any $r\ge 0$, the $n\times n$ matrix
$(f^r(x_i, x_j))$ is positive semi-definite. So the matrix $(f(x_i,
x_j))$ is infinitely divisible. One can easily check the following
identities:
$$
\bigg(\frac{1}{f^r[x_i, x_j]}\bigg)=D(f^r(x_i, x_j))D, \
\bigg(\frac{f^r(x_i, x_j)}{f^r[x_i, x_j]}\bigg)=D(f^{2r}(x_i,
x_j))D,
$$
where $D={\rm diag}(\frac{1}{f^r(x_1)}, ..., \frac{1}{f^r(x_n)})$.
Thus $(\frac{1}{f^r[x_i, x_j]})$ and $(\frac{f^r(x_i, x_j)}{f^r[x_i,
x_j]})$ are positive semi-definite if $r\ge 0$. In other words, the
matrices $(\frac{1}{f[x_i, x_j]})$ and $(\frac{f(x_i, x_j)}{f[x_i,
x_j]})$ are infinitely divisible. \\

{\sc Remark.} It should be pointed out that the condition that $f$
is multiplicative is necessary. Otherwise the conclusion may be
false. For example, let $S=\{6, 10, 15\}$ and $f$ be defined by
$f(1)=f(3)=0, f(10)=3$ and $f(m)=1$ for $m\ne 1, 3, 10$. It is clear
that $f\in {\cal C}_S$ and $f$ is not multiplicative. Obviously,
$(f(x_i, x_j))=\left(\begin{array}{ccc}
  1 & 1& 0 \\
  1 & 3& 1 \\
  0 & 1 & 1
\end{array}\right)$ is positive semi-definite. One can easily prove that $(f(x_i,
x_j))^{\circ r}$ is positive semi-definite if and only if $r\ge \log
2/\log 3$. Hence the matrix $(f(x_i, x_j))$ is not infinitely
divisible.\\

Let $l\ge 0$ be an integer. For any arithmetical function $f$,
define its $l$-th Dirichlet convolution, denoted by $f^{(l)}$,
inductively as follows: $f^{(0)}:=\delta$ and $f^{(l)}:=f^{(l-1)}*f$
if $l\ge 1$, where $\delta$ is the function defined for any positive
integer $m$ by $\delta(m):=1, \ {\rm if} \  m=1; 0, \ {\rm
otherwise}.$ Evidently, $f*\delta =f$ for any arithmetical function
$f$ and $f^{(l)}:=\underbrace{f*...*f}_{l\ {\rm times}}.$ Let ${\bf
Z}_{>0}$ denote the set of positive integers.\\

{\sc Lemma 2.3.} {\it Let $c\ge 1$ and $d\ge 0$ be integers. If
$f_1, ..., f_c\in {\cal C}_{S}$ are distinct arithmetical functions
and $(l_1,...,l_c)\in {\bf Z}_{>0}^c$ satisfies $l_1+...+l_c>d$,
then $f_1^{(l_1)}*...*f_c^{(l_c)}*\mu ^{(d)}\in {\cal C}_{S}$.}

{\it Proof.} Clearly to prove Lemma 2.3 it is sufficient to prove
that for any integer $l>d$ and any (not necessarily distinct)
arithmetical functions $g_1, ..., g_l\in {\cal C}_{S}$, we have $
g_1*...*g_l*\mu ^{(d)}\in {\cal C}_{S}.$ In the following let $g_1,
..., g_l\in {\cal C}_{S}$ and $l>d$. Now for any $x\in S$ and any
$m|x$, since $l\ge d+1$, we have
$$
\begin{array}{rl}
&((g_1*...*g_l*\mu ^{(d)})*\mu )(m)\\
=&(g_1*...*g_l*\mu^{(d+1)})(m)\\
=&((g_1*\mu)*...(g_d*\mu )*(g_{d+1}*\mu )*g_{d+2}*...*g_l)(m)\\
=&\displaystyle \sum_{\SSS{m_1...m_l=m}{(m_1,...,m_l)\in {\bf
Z}_{>0}^{l}}}(g_1*\mu)(m_1)...(g_{d+1}*\mu)(m_{d+1})g_{d+2}(m_{d+2})...g_l(m_l).
\end{array} \eqno(5)
$$

For $1\le i\le d+1$, since $g_i\in {\cal C}_{S}$ and $m_i|x$, we
have $(g_i*\mu)(m_i)\ge 0.$ On the other hand, for $d+2\le j\le l$,
$g_j\in {\cal C}_{S}$ together with $m_j|x$ implies that
$g_j(m_j)=\sum_{d'|m_j}(g_j*\mu )(d')\ge 0.$ It then follows from
(5) that $ ((g_1*...*g_l*\mu^{(d)})*\mu )(m)\ge 0$ as desired.\\

From Theorem 2.1 and Lemma 2.3 we deduce immediately that the
following more general result is true.\\

{\sc Theorem 2.2.} {\it Let $c\ge 1$ and $d\ge 0$ be integers. If
$f_1, ..., f_c\in {\cal C}_{S}$ are distinct and multiplicative and
$(l_1,...,l_c)\in {\bf Z}_{>0}^c$ satisfies $l_1+...+l_c>d$, then
the following three $n\times n$ matrices
$$((f_1^{(l_1)}*...*f_c^{(l_c)}*\mu
^{(d)})(x_i, x_j)), \
\bigg(\frac{1}{(f_1^{(l_1)}*...*f_c^{(l_c)}*\mu ^{(d)})[x_i,
x_j]}\bigg),$$
$$\bigg(\frac{(f_1^{(l_1)}*...*f_c^{(l_c)}*\mu
^{(d)})(x_i, x_j)}{(f_1^{(l_1)}*...*f_c^{(l_c)}*\mu ^{(d)})[x_i,
x_j]}\bigg)$$
are infinitely divisible.}\\

In particular, we have\\

{\sc Theorem 2.3.} {\it If $f\in {\cal C}_{S}$ is multiplicative and
$l>d\ge 0$ are integers, then the following three $n\times n$
matrices
$$((f^{(l)}*\mu
^{(d)})(x_i, x_j)), \ \bigg(\frac{1}{(f^{(l)}*\mu ^{(d)})[x_i,
x_j]}\bigg), \ \bigg(\frac{(f^{(l)}*\mu ^{(d)})(x_i,
x_j)}{(f^{(l)}*\mu ^{(d)})[x_i, x_j]}\bigg)$$
are infinitely
divisible.}\\

{\bf 3. Examples.} In this section, we give some examples to
illustrate our main results.\\

{\sc Example 3.1.} Let $\xi _\varepsilon$ be defined by $\xi
_\varepsilon(m)=m^\varepsilon$ for any integer $m\ge 1$. It is easy
to check that $\xi _\varepsilon\in {\cal C}_S$ for any set $S$ of
positive integers and any $\varepsilon\ge 0$. By Theorem 2.1, the
matrices
$$
((x_i, x_j)^\varepsilon), \ \bigg(\frac{1}{[x_i,
x_j]^\varepsilon}\bigg), \ \bigg(\frac{(x_i, x_j)^\varepsilon}{[x_i,
x_j]^\varepsilon}\bigg)
$$
are infinitely divisible for any nonnegative real numbers
$\varepsilon$. Note that Hong-Loewy [14], Hong-Lee [13] and
Lindqvist and Seip [16] investigated the asymptotic behavior of the
eigenvalues of the above three matrices respectively. \\

{\sc Example 3.2.} Let $J_\varepsilon:=\xi_\varepsilon *\mu $ be the
generalized Jordan function. Since $(J_{\varepsilon}*\mu
)(p)=p^\varepsilon -2\ge 2^\varepsilon -2\ge 0$ for any prime $p$
and any real number $\varepsilon\ge 1$, we have $J_{\varepsilon}\in
{\cal C}_S$ for any set $S$ of positive integers and any
$\varepsilon\ge 1$. By Theorem 2.1, the matrices
$$
(J_{\varepsilon}(x_i, x_j)), \ \bigg(\frac{1}{J_{\varepsilon}[x_i,
x_j]}\bigg), \ \bigg(\frac{J_{\varepsilon}(x_i,
x_j)}{J_{\varepsilon}[x_i, x_j]}\bigg)
$$
are infinitely divisible for any real number
$\varepsilon \ge 1$.\\

{\sc Example 3.3.} Let $\xi _\varepsilon$ and
$J_\varepsilon:=\xi_\varepsilon *\mu $ be defined as above. Then by
Theorem 2.2, the matrices
$$
((\xi _\varepsilon^{(l)}*J_{\epsilon}^{(t)}*\mu ^{(d)})(x_i, x_j)),
\ \bigg(\frac{1}{(\xi _\varepsilon^{(l)}*J_{\epsilon}^{(t)}*\mu
^{(d)})[x_i, x_j]}\bigg), \ \bigg(\frac{(\xi
_\varepsilon^{(l)}*J_{\epsilon}^{(t)}*\mu ^{(d)})(x_i, x_j)}{(\xi
_\varepsilon^{(l)}*J_{\epsilon}^{(t)}*\mu ^{(d)})[x_i, x_j]}\bigg)
$$
are infinitely divisible for any real numbers $\varepsilon\ge 0$ and
$\epsilon\ge 1$ and any nonnegative integers $l, t$ and
$d$ such that $l+t>d$.\\

Finally, we remark that Hong [12] and Hong and Loewy [15] studied
the asymptotic behavior of the eigenvalues of matrices associated
with arithmetical functions including all the matrices in Examples
3.1-3.3 as special examples.\\

{\bf Acknowledgement.} The author would like to thank the anonymous
referee for valuable comments and suggestions.

\noindent{Mathematical College}\\
Sichuan University\\
Chengdu 610064, P.R. China\\
E-mail: sfhong@scu.edu.cn\\
s-f.hong@tom.com\\
hongsf02@yahoo.com

\end{document}